# A Novel Mathematical Formulation for Density-based Topology Optimization Method Considering Multi-Axis Machining Constraint

Hao Deng, Albert C. To[*]

Department of Mechanical Engineering and Materials Science, University of Pittsburgh, Pittsburgh, PA 15261

*Corresponding author. Email: albertto@pitt.edu

**Abstract**

This paper proposes a novel density-based method for structural design considering restrictions of multi-axis machining processes. A new mathematical formulation based on Heaviside function is presented to transform the design field into a geometry which can be manufactured by multi-axis machining process. The formulation is developed for 5-axis machining, which can be also applied to 2.5D milling restriction. The filter techniques are incorporated to effectively control the minimum size of void region. The proposed method is demonstrated by solving the compliance minimization problem for different machinable freeform designs. The length to diameter (L:D) ratio geometric constraint is introduced to ensure the machinable design, where deep hole or narrow chamber features are avoided using proposed method. Several two and three-dimensional numerical examples are presented and discussed in detail.

## 1. Introduction

Topology optimization has become an important tool to determine the optimal shape for maximum performance subject to given design constraints. These performance objective functions and design constraints affect the applicability of the design and typically include the mass, compliance, stress, or displacement due to loading conditions. Research for modern topology optimization first started with the work of Bendsoe and Kikuchi [1], and new methods have since been implemented such as the Solid Isotropic Material with Penalization (SIMP) method [2] and the level-set method [3], both of which can be found in commercial software. For the SIMP method, each element is assigned a density variable to control the material distribution and remove material from gradient-based optimization. For the level-set method, the boundary of the shape geometry is the zero level-set of an implicit function. The parameters defining the function are optimized to produce the optimal shape.

These methods can produce complex geometries not typically considered by a human engineer or designer. Consequently, it can be difficult to manufacture these designs and may require extensive post-processing [4]. To control the geometric complexity, other new methods have been developed such as moving morphable components (MMC) [5, 6], moving morphable voids (MMV) [7], and geometry projection [8-10]. MMC and MMV use controlled geometries to add or remove material from the design domain and determine the optimal shape. Geometry projection methods [10, 11] use geometries to describe the density distribution on a background mesh. Additionally, integrating the topology optimization process with CAD has been addressed as well to further reduce post-processing [12, 13].

This paper mainly focuses on the multi-axis machining constraints, which is a widely used techniques in subtractive manufacturing for metal component production. For multi-axis machining, the relative position and orientation of cutting tools and workpiece can be manipulated in 4 or 5 degrees of freedom. The unnecessary material is removed by machining tool until the desired shape achieves. Compared with 2.5D milling, the multi-axis machining allows more design freedom. Besides machining, metal additive manufacturing is also a well-known technique to print free-form design with high design freedom, while the material strength and fatigue properties of parts made of AM still have large gap compared to metal machining. For high-strength aerospace or naval structures, multi-axis machining is an ideal option for manufacturing parts. Thus, topology optimization approaches considering machining constraint is a necessary and valuable research to achieve a trade-off between the available manufacturing technologies. In recent years, a few effective approaches for multi-axis machining-based topology optimization are proposed. For 2.5D milling, Liu et al [14] proposes an explicit feature-based level-set method, where the feature fitting algorithm is incorporated into the boundary evolvement process. Furthermore, Liu and Albert [15] proposed a novel CAD-based topology optimization system for milling constraint, where feature and dynamic modeling history is incorporated into the optimization process. For multi-axis machining, Amir et al [16] presented a topology optimization framework using convolutions in configuration space to enable manufactured design using multi-axis machining, where an inaccessibility measure field of design domain is introduced to identify non-manufacturable features. Nigel et al [17] proposed a level-set-based topology optimization method for multi-axis machining, where the advection velocity at every iteration is modified to ensure the manufacturability conditions. This level-set-based method simulates the subtractive process by cutting materials accessible to the machine in every iteration until optimization converges. Matthijs et al [18] proposed a density-based approach incorporating multi-axis restrictions, where a filter technique based on KS aggregation function [19] is introduced to transform an input design field into a manufacturable geometry. Amir et al [20] proposed framework to optimize the build orientation with respect to removability of support structures, where a Pareto-optimality criterion is implemented to achieve a trade-off between build orientation and machining setup. Meanwhile, Amir et al [21] further extended the notion of inaccessibility measure field (IMF) to support's structures of additive manufacturing to identify the inaccessible points, and the IMF is integrated into the sensitivity field to guide the TO results. This approach enables optimal design can be additive manufactured and post-processed for support removal using machining techniques. Besides the multi-axis machining, several other related topology optimization studies regarding casting constraints can also be found in Refs. [22-25].

The aim of this paper is to propose a novel density-based approach to optimize parts for multi-axis machining restrictions, where a concise aggregation-free density projection formulation is demonstrated. The proposed method has concise mathematical expression without aggregation function, which is easy for sensitivity analysis and numerical implementation. The length to diameter (L:D) ratio constraint is taken into consideration in the proposed scheme. The current method is in the framework of density-based method, which is able to inherit the existing sensitivity from standard density method with chain rule for different physical problems. The subtractive manufacturing constraints can be satisfied after density projection. The remainder of this paper is organized as follows. Section 2 presents the mathematical formulation of density projection algorithm. Section 3 demonstrates three typical topology optimization problems with machining constraint from 2D to 3D to verify the effectiveness of the proposed method, followed by conclusions in Section 4.

## 2. Mathematical Formulation for Multi-axis Machining

## 2.1 Machining Restriction

The milling head is generally described by a cylinder capped with a hemisphere oriented in the milling direction. The end of the head cylinder is assumed to extend infinitely far away from the cutting surface. For 3-axis machining, the workpiece is still while the cutting tool moves along the 3 mutually perpendicular axes to mill the part, which is the most widely used technique to create a metal part. While for 5-axis machining process, the milling tool can be oriented to approach the surface from an arbitrary direction. The machinable part is determined by whether all its surface points are accessible by a tool bit without any intersection with the interior of the part. As described by Ref. [18], machinability means that the density field must be monotonically increasing in the insertion direction without any internal holes. To enforce the milling constraints, a Heaviside function-based projection method is proposed here to construct the monotonically increasing density field in the milling direction.

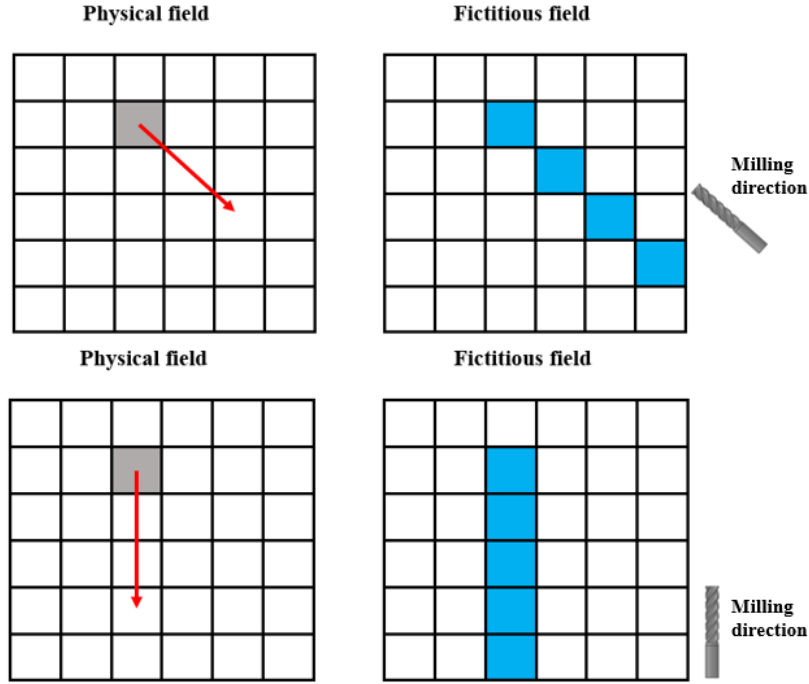

Figure 1. Density Mapping from fictitious field to physical field

This subsection outlines the conceptual idea for the proposed mathematical formulation considering multi-axis machining constraint. The physical field is used for finite element analysis. The physical density field is computed through the cumulative summation of fictitious density along reverse milling direction as shown in Fig. 1. The Heaviside function is introduced here to map the fictitious density to physical density as follows,

$$\rho_p = H_2\left(\Sigma_{j \in \mathcal{M}} H_1(\rho_f^{(j)})\right) \qquad (1)$$

where notation $\Sigma(\cdot)$ means cumulative summation of fictitious element density along reverse milling direction. $\mathcal{M}$ is the set of elements along reverse milling direction. The element set $\mathcal{M}$ can be determined based on following flowchart,

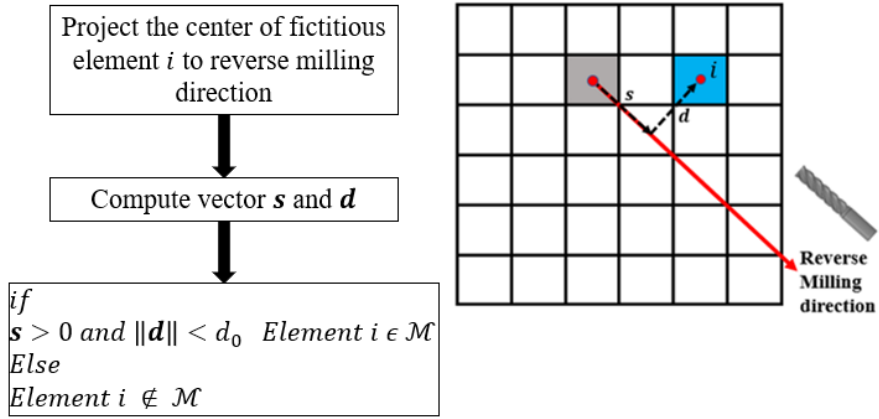

Figure 2. Flowchart of computing set $\mathcal{M}$

As shown in flowchart, $d_0$ is a threshold, and the operation $\|\cdot\|$ denotes the Euclidean norm. Note that the computational cost will increase with the increase of the element number. A table is listed below to demonstrate the computational costs for rectangle domain based on MATLAB code [16] for one milling direction. The element set $\mathcal{M}$ is just computed once and store the information in a matrix, which will be reused during optimization.

Table 1. Computational Cost for Element set $\mathcal{M}$ (single machining direction [1,1,1])

| Element number | 10 × 10 × 10 | 25 × 25 × 25 | 50 × 50 × 50 | 75 × 75 × 75 | 100 × 100 × 100 |
|---|---|---|---|---|---|
| Time cost | 0.013s | 0.21s | 8.53s | 18.78s | 89.52s |

As demonstrated in Fig. 1, the grey element in physical field $\rho_p$ is computed through summation of blue element ($\mathcal{M}$) in fictitious density field $\rho_f$ according to Eq. (1). $H_1(\cdot)$ and $H_2(\cdot)$ denote the approximate Heaviside functions defined as:

$$\begin{cases} H_1(x) = \frac{\tanh(10x-3)+1}{2} \\ H_2(x) = \frac{\tanh(6x-3)+1}{2} \end{cases} \quad (2)$$

The shape of Heaviside functions is plotted as follows,

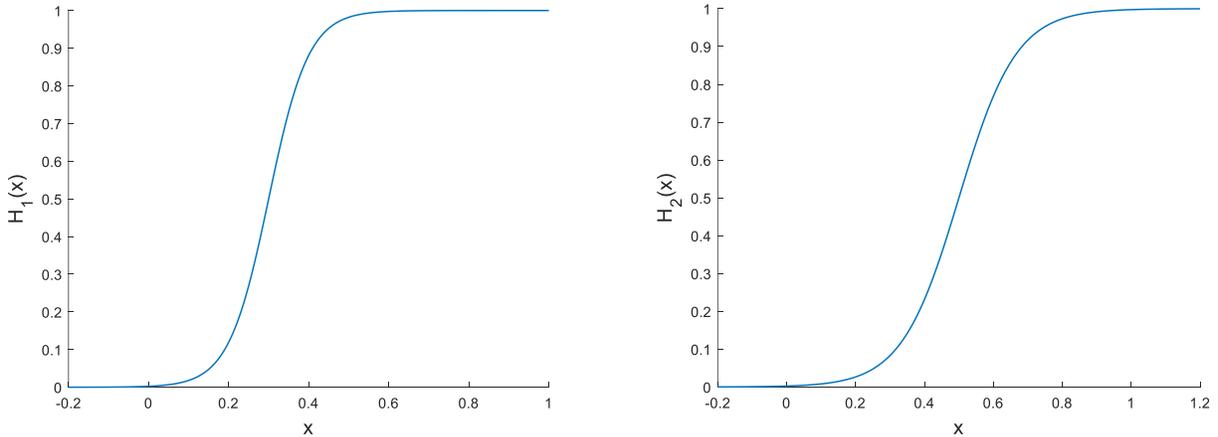

Figure 3. Heaviside functions

The reason for introducing two Heaviside functions with different parameters is due to numerical stability during optimization process after numerical experiments. In general, low coefficient values are preferred, which will make the objective function smooth enough. However, a low coefficient value will weaken the filter effect, which is critical to enable the machinable design. Meanwhile, a higher coefficient value will generate less intermediate density in final optimal designs. Thus, a reasonable coefficient for Heaviside function need to be chosen to achieve a trade-off between smoothness and filter effects.

To effectively control the minimum size of void area, we introduce the void field to work as design variables. The relationship between void field $\rho_v$ and fictitious domain $\rho_f$ can be defined as follows:

$$\rho_f = \frac{1}{\sum_{i \in N_f} d_{fi}} \sum_{i \in N_f} d_{fi} \left(1 - \rho_v^{(i)}\right) \tag{3}$$

Note that $\rho_v = 1$ represent the void element, while $\rho_v = 0$ means solid element. $(\cdot)^{(i)}$ means the $i_{th}$ element. $d_{fi}$ is weight factor for filter defined as:

$$d_{fi} = max(0, r_{min} - \Delta(f, i)) \tag{4}$$

$r_{min}$ is the filter radius. $\Delta(f, i)$ denotes the distance between the center of element $i$ and $f$. $N_f$ is the set of elements $i$ for which the distance $\Delta(f, i)$ is less than $r_{min}$. More details regarding filtering techniques can be found in Ref. [26]. Therefore, the physical density $\rho_p$ can be explicitly expressed as,

$$\rho_p = H_2 \left( \sum_{f \in \mathcal{M}} H_1 \left( \frac{1}{\sum_{i \in N_f} d_{fi}} \sum_{i \in N_f} d_{fi} \left(1 - \rho_v^{(i)}\right) \right) \right) \tag{5}$$

For multiple milling directions ($k = 1, \cdots, n$), the composite physical density field can be simply expressed as,

$$\overline{\rho_p} = \prod_{k=1}^{n} \left(\rho_p^{(k)}\right) \tag{6}$$

where $\rho_p^{(k)}$ denote the physical density computed from the $k$ milling direction. A typical projection from void density field to physical density field is demonstrated in Fig. 4.

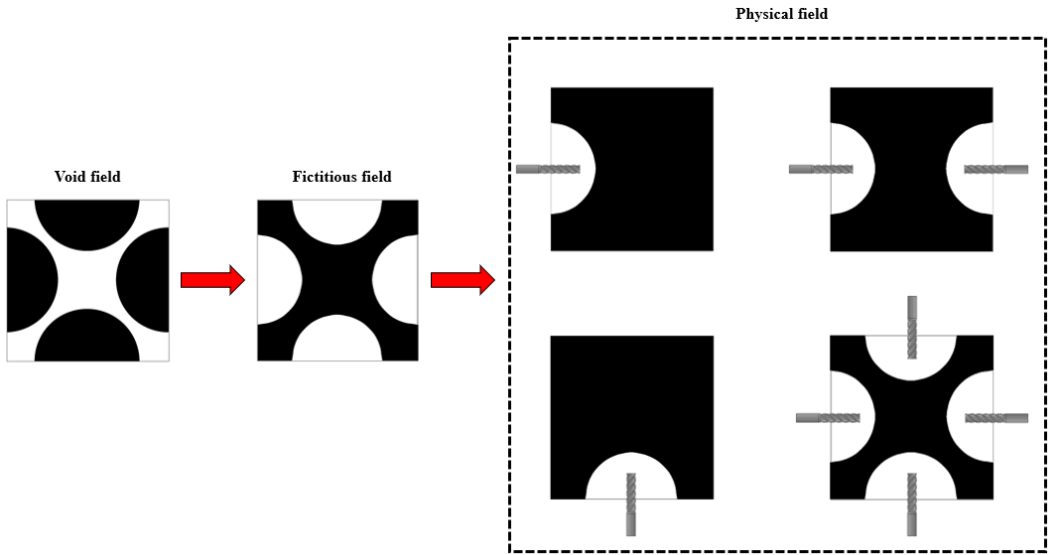

Figure 4. Projection from the void field to physical field

## 2.2 Sensitivity Analysis

Given a response $f(\rho_v)$, the sensitivity with respect to the design variables $\rho_v$ can be obtained through the chain rule as follows,

$$\frac{\partial f}{\partial \rho_v^{(l)}} = \frac{\partial f}{\partial \overline{\rho_p}} \frac{\partial \overline{\rho_p}}{\partial \rho_v^{(l)}} = \frac{\partial f}{\partial \overline{\rho_p}} \sum_{k=1}^{n} \left( \frac{\partial \rho_p^{(k)}}{\partial \rho_v^{(l)}} \prod_{j=1; j \neq k}^{n} \left( \rho_p^{(j)} \right) \right) \tag{7}$$

where Einstein summation applies to all repeated subscripts. $n$ is the number of machining direction. The explicit form of the term $\frac{\partial \rho_p^{(i)}}{\partial \rho_v^{(l)}}$ can be expressed as,

$$\frac{\partial \rho_p^{(i)}}{\partial \rho_v^{(l)}} = H_2' \cdot \left( \sum_{f \in \mathcal{M}^{(i)}} \frac{\partial H_1 \left( \frac{1}{\sum_{i \in N_f} d_{fi}} \sum_{i \in N_f} d_{fi} \left(1 - \rho_v^{(i)}\right) \right)}{\partial \rho_v^{(l)}} \right) \tag{8}$$

where

$$\frac{\partial H_1 \left( \frac{1}{\sum_{i \in N_f} d_{fi}} \sum_{i \in N_f} d_{fi} \left(1 - \rho_v^{(i)}\right) \right)}{\partial \rho_v^{(l)}} = H_1' \cdot \frac{\partial \left( \frac{1}{\sum_{i \in N_f} d_{fi}} \sum_{i \in N_f} d_{fi} \left(1 - \rho_v^{(i)}\right) \right)}{\partial \rho_v^{(l)}} \tag{9}$$

where $H_1'$ and $H_2'$ denote the first derivative of Heaviside function defined in Eq. (2). The detailed form of $\frac{\partial f}{\partial \rho_v^{(l)}}$ can be easily derived, which is omitted here due to space limitations.

## 2.3 Machining Restriction considering length to diameter (L:D) ratio

Besides the density field must be monotonically increasing in the insertion direction. The other crucial machinable constraint is the length to diameter (L:D) ratio [27]. The (L:D) ratio is an important concept for machinable design, which has not been considered in the previous topology optimization research. As shown in Fig. 5, the hole is too deep for the machine tool to reach. Although this case satisfies the density non-decreasing requirement along the insertion direction, the part still cannot be machinable due to the inaccessible deep hole. This situation can be resolved if the cone angle geometrical constraints are introduced. As shown in Fig. 6, a cylinder cone angle $\theta$ is defined to further constrain design as follows,

$$\theta = actan\left(\frac{D}{L}\right) \quad \left[0, \frac{\pi}{2}\right] \tag{10}$$

where $D$ is the diameter and $L$ is the length as shown in Fig. 6. This geometrical constraint means that there should be no material inside the circular cone envelop surface (CES) measured from arbitrary point on the design external surface. This CES geometry constraint can effectively avoid feature with small (L:D) ratio, e.g. deep holes.

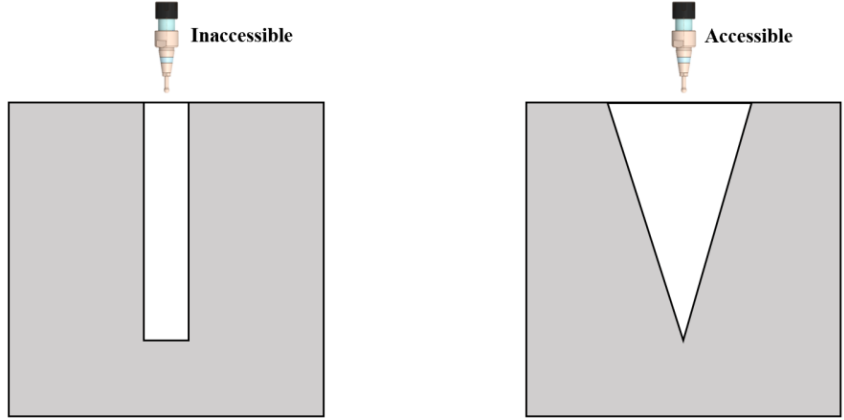

Figure 5. Machining Restriction considering length to diameter (L:D) ratio

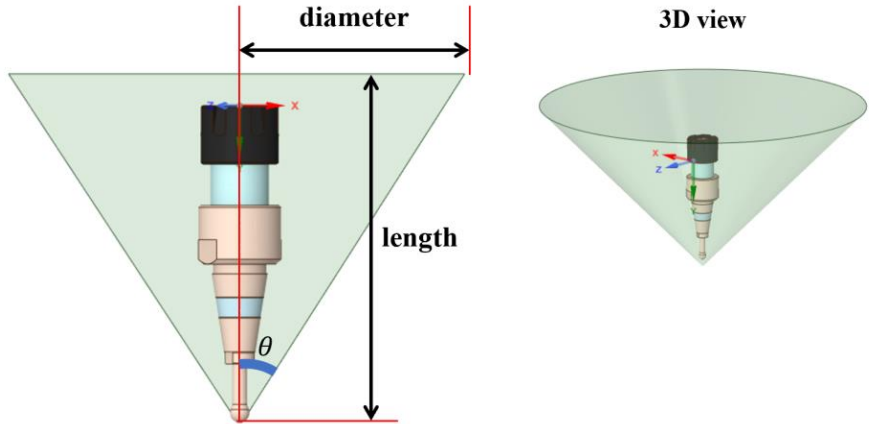

Figure 6. Machining Restriction considering length to diameter (L:D) ratio

To apply the CES geometric constraint in topology optimization, the mathematical formulation for density mapping (Eq. (5)) need to be modified as follows,

$$\rho_p = H_2 \left( \sum_{f \in \mathcal{M}_{CES}} H_1 \left( \frac{1}{\sum_{i \in N_f} d_{fi}} \sum_{i \in N_f} d_{fi} \left(1 - \rho_v^{(i)}\right) \right) \right) \quad (11)$$

The element set $\mathcal{M}_{CES}$ can be determined based on following equation,

$$\arctan\left(\frac{\|\mathbf{d}\|}{h}\right) \leq \theta \quad (12)$$

where $\theta$ is CES angle, which will be determined by CNC machine and tool shape. $\|\mathbf{d}\|$ is the minimum distance from centroid of the element to red line as shown in Fig. 7, and $h$ is the projection length. The sensitivity analysis process is similar to section 2.2, which is omitted here.

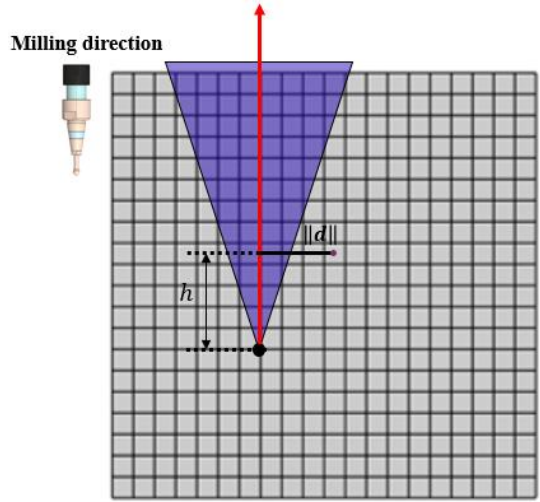

Figure 7. Density Mapping considering CES constraint

A simple schematic is shown here to demonstrate the difference between two projection schemes. As shown in Fig. 8, the projection design mapped from fictitious field generates small deep hole if no CES constraint is applied, while no deep hole appears once the CES constraint ($\theta = 30°$) is implemented.

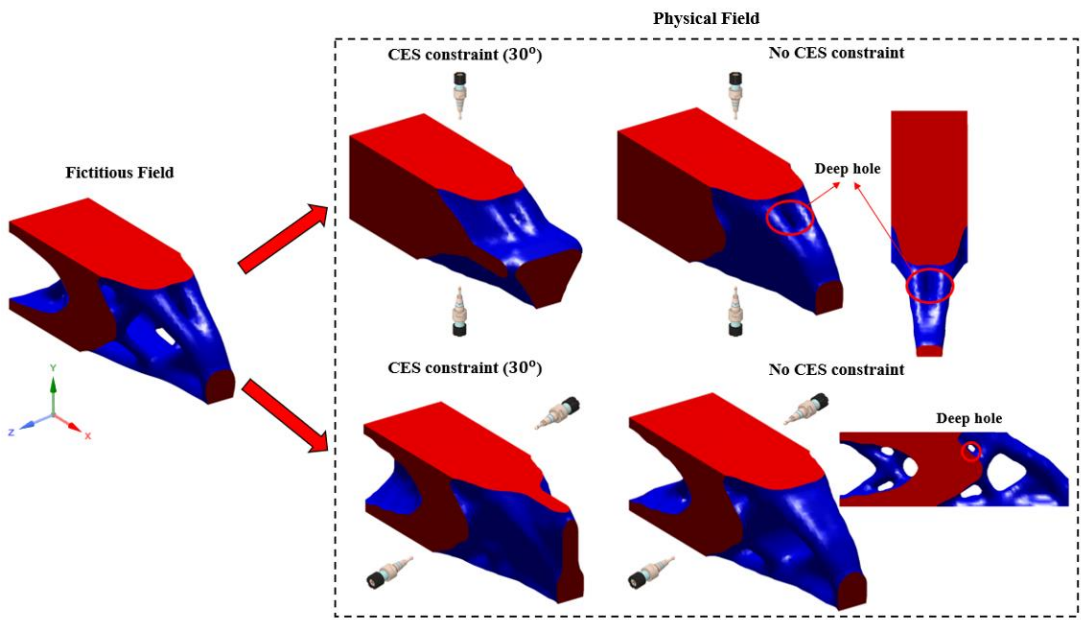

Figure 8. Physical field computed by two projection schemes

## 2.4 Objective and Constraints

In this section, the compliance minimization optimization problem can be formulated as following,

$$\begin{cases} \text{Find:} \rho_v \\ \text{Min:} C(u, \rho_v) = \frac{1}{2} \int_{\Omega} \boldsymbol{\varepsilon}(\boldsymbol{u})^T \boldsymbol{D}(\rho_p(\rho_v))\boldsymbol{\varepsilon}(\boldsymbol{u}) d\Omega \\ S.t: \left\{ \frac{1}{|\Omega|} \int_{\Omega} \rho_p(\rho_v) d\Omega - V_{prescribe} \leq 0 \right. \end{cases}$$

where $\rho_v$ is the design variables, and $C$ is the objective defined by the structural compliance. $\rho_p$ is the physical density field, $\Omega$ the design domain, and $V_{prescribe}$ the prescribed volume fraction. $\boldsymbol{u}$ is the unknown displacement field and $\boldsymbol{\varepsilon}$ is the strain resolved by FEM, and $\boldsymbol{D}$ is the elastic tensor matrix. The sensitivity of objective function and constraint with respect to design variable $\rho_v$ can be obtained through the chain rule.

## 3. Numerical Examples

### 3.1 Machining-based Optimization for a 2D cantilever beam

The first 2D numerical example for machining-constrained topology optimization is demonstrated in Fig. 5. The design domain is uniformly meshed by $100 \times 100$ quad elements with unit length. The loading $F = 1$ is applied on the right-bottom corner, and left side is fully fixed. The volume fraction constraint $\bar{V}$ is chosen as 0.2. The elastic constants are chosen as follows: Elastic modulus $E = 1$ and Poisson's ratio $\mu = 0.3$. The filter radius for design variable is chosen as $r_{min} = 4$. The initial void field (design variable) is plotted in Fig. 3. To make a comparison, different milling direction constraints are applied to produce an optimized design. The reference solution without any manufacturing constraints is demonstrated in Fig. 6. Obviously, the reference solution cannot be manufactured if machining operations are limited in the $x, y$-plane. The initial and optimized designs for single milling orientation are demonstrated in Fig. 7, and convergence history is plotted in Fig. 8. Compared with the reference solution, the optimized part for unidirectional milling restriction shows strong limitation of design freedom, which also has great impact on structural compliance. Note that the designs with single milling orientation constraint are considerably more compliant with respect to reference design. For multiple tool orientations, the optimized designs are shown in Fig. 9, and convergence history is plotted in Fig. 10. The structural performance with multiple directions is better than the single orientation, while the compliance value is still higher than the reference design. This simple 2D numerical example proves the effectiveness of proposed method to force the optimization process towards a different solution considering multi-axis machining constraints.

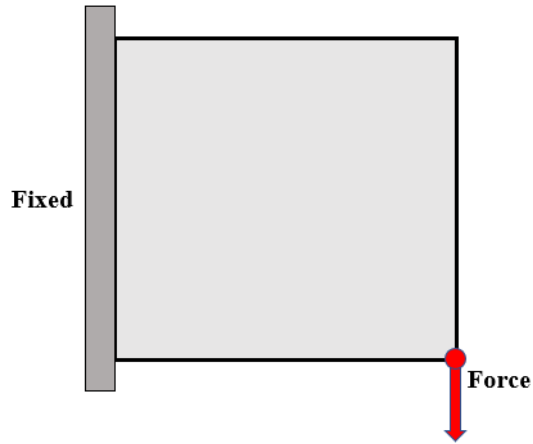

Figure 5. Two-dimensional Cantilever Beam

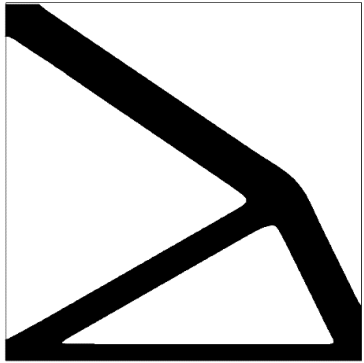

Figure 6. Reference design (Compliance: 50.05)

(a) Initial design                                   (b) Optimized design

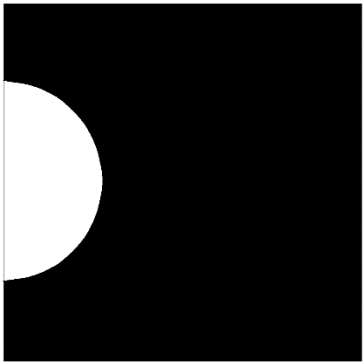
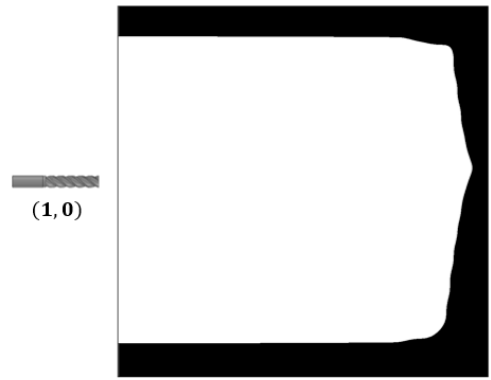

(a) Initial design                                   (b) Optimized design

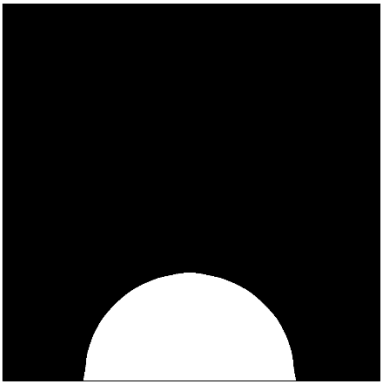
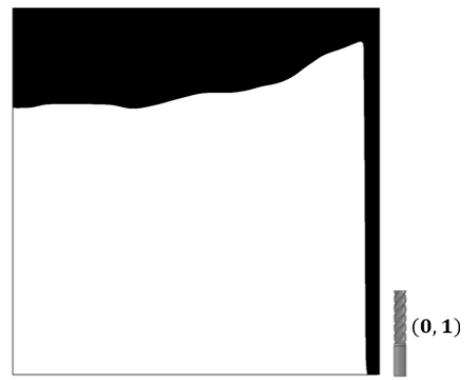

Figure 7. Designs obtained using single tool orientation

(a) Orientation (1,0)                              (b) Orientation (0,1)

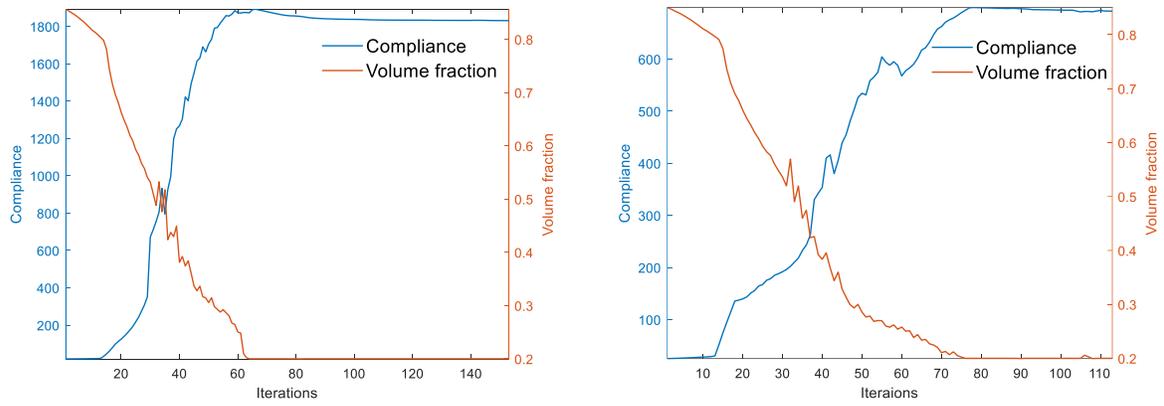

Figure 8. Convergence history

(a) Initial design            (b) Optimized design

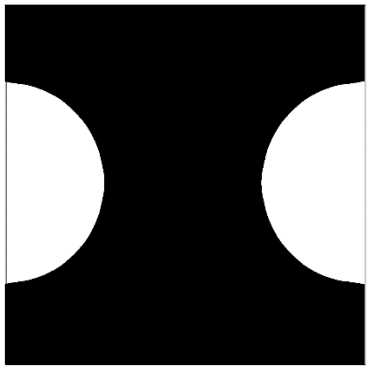 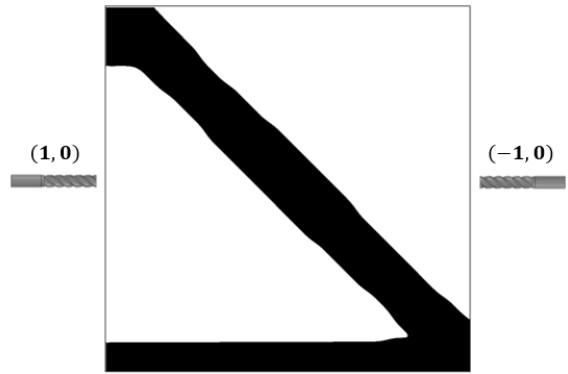

(a) Initial design            (b) Optimized design

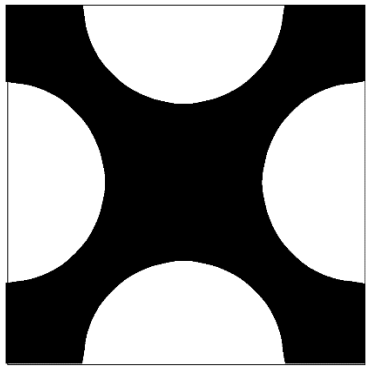 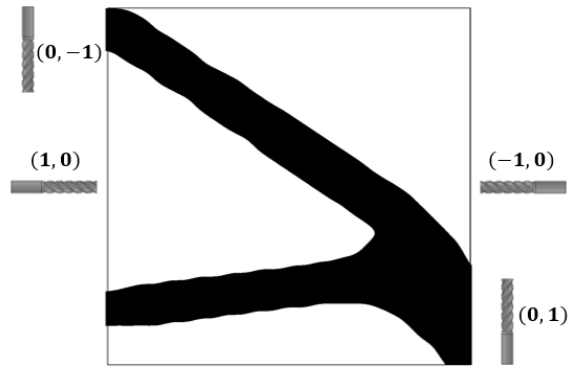

Figure 9. Designs obtained using multiple tool orientations

(a) 2 Orientations            (b) 4 Orientations

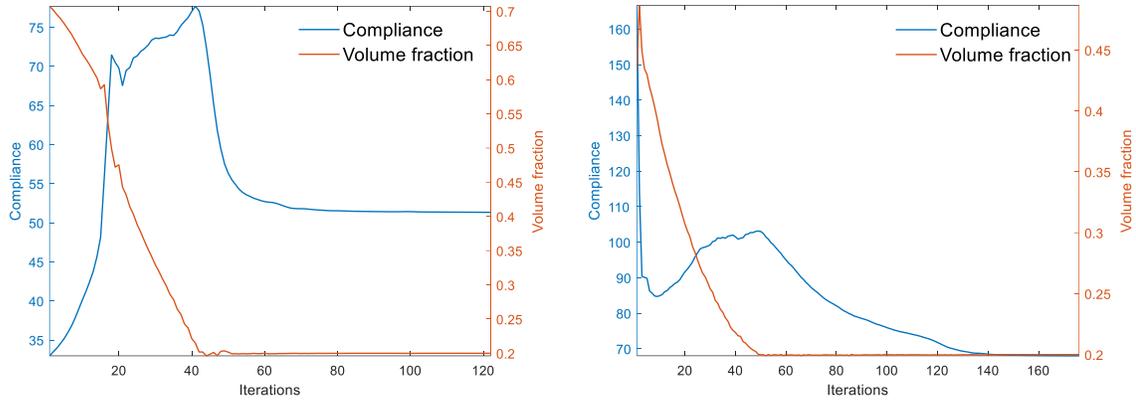

Figure 10. Convergence history

### 3.2 Machining-based Optimization for 3D Cantilever beam

#### 3.2.1 Optimized Design without CES constraint

In this subsection, we focus on verifying the effectiveness of proposed method for 3D designs. In the first test example, a three-dimensional cantilever beam example is presented for compliance minimization. The unit force is applied at the midpoint of right-bottom edge, and the left end is fully fixed. The design domain is discretized by a $144 \times 48 \times 48$ hexahedral mesh with unit element length. The elastic constants are chosen as follows: Elastic modulus $E = 1$ and Poisson's ratio $\mu = 0.3$. The volume fraction constraint is set to be 0.3. The filter radius for multi-axis machining optimization is chosen as $r_{min} = 4$. A reference solution is demonstrated in Fig. 13, where no machining constraints are applied. As shown in Fig. 13, the reference optimized design consists of hollow chamber, which is not accessible by milling tools. It is worth to mention that the initialization of fictitious field for optimization is shown in Fig. 12. To ensure manufacturability of the optimized design, the results of two-orientation machining optimization for the cantilever beam are plotted in Fig. 8. The orientation of machine tool is described by the vector $(a, b, c)$ in the local coordinate as shown in Fig. 11. The compliance value of designs obtained from different orientations (Fig. 14) is close and slightly higher than the reference design. For multi-axis machining constraints, the optimized result is demonstrated in Fig. 15 with 6 different milling orientations. However, compliances differ only slightly among all designs, which denotes the number of different directions has little effect on the current issue.

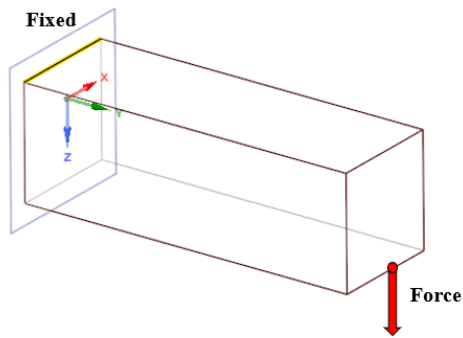

Figure 11. Three-dimensional Cantilever Beam

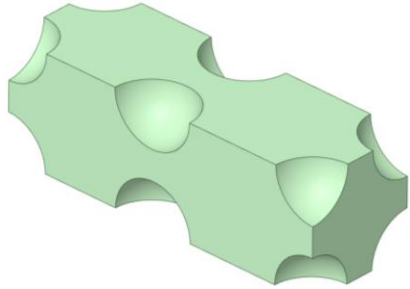

Figure 12. Initialization of fictitious field

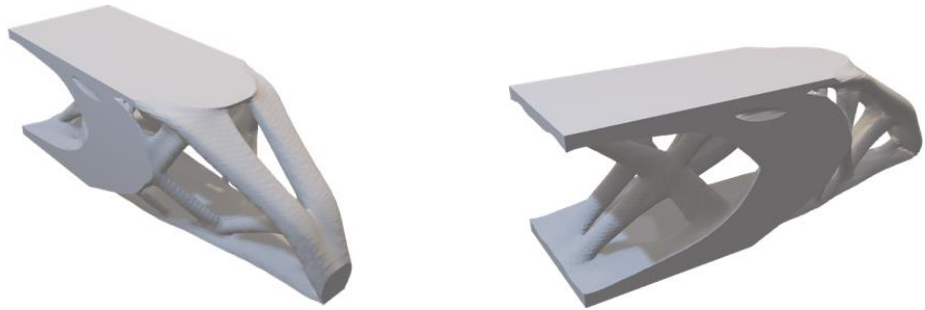

Figure 13. Reference solution without machining constraints (Compliance:10.43)

(a) Compliance:11.18

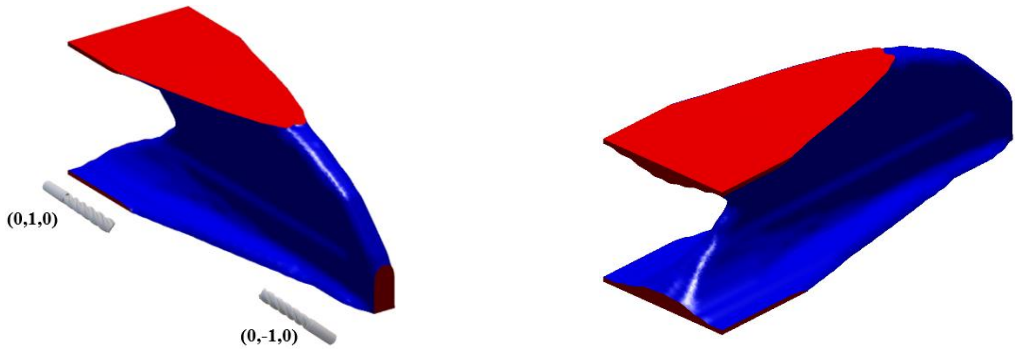

(b) Compliance:10.77

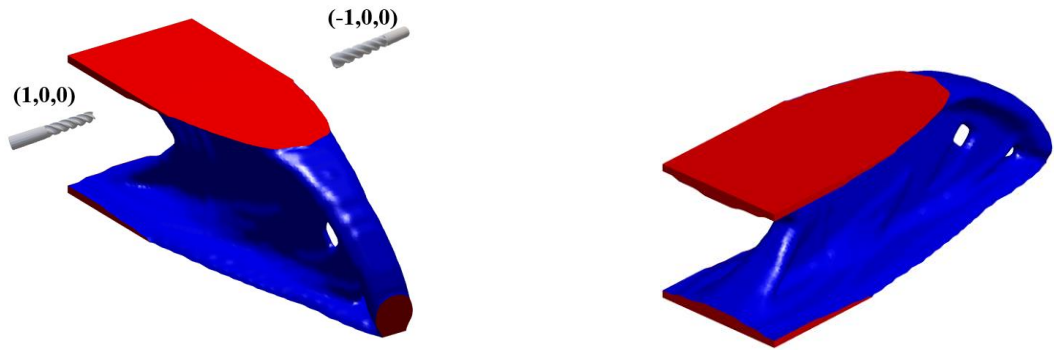

(c) Compliance:13.64

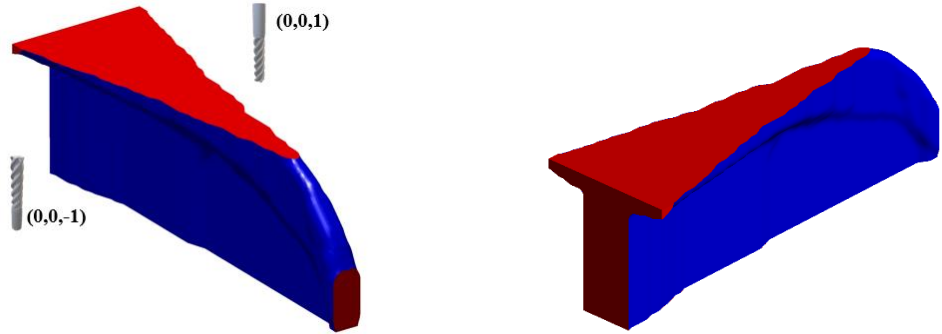

(d) Compliance:12.04

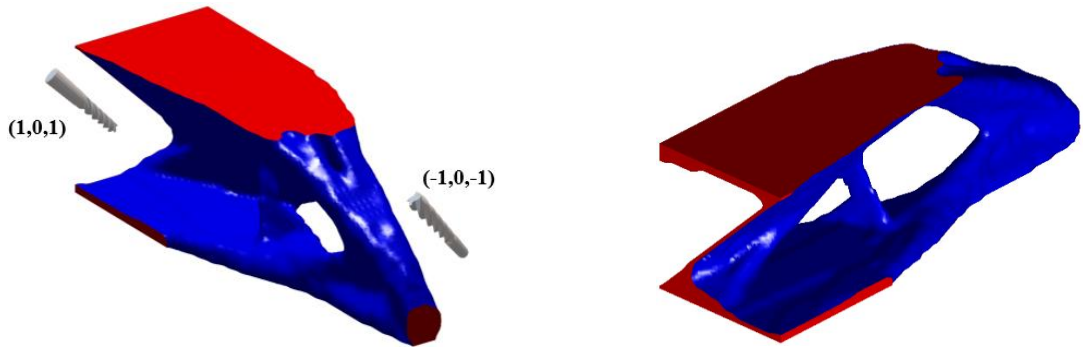

Figure 14. Two orientation machining constraints

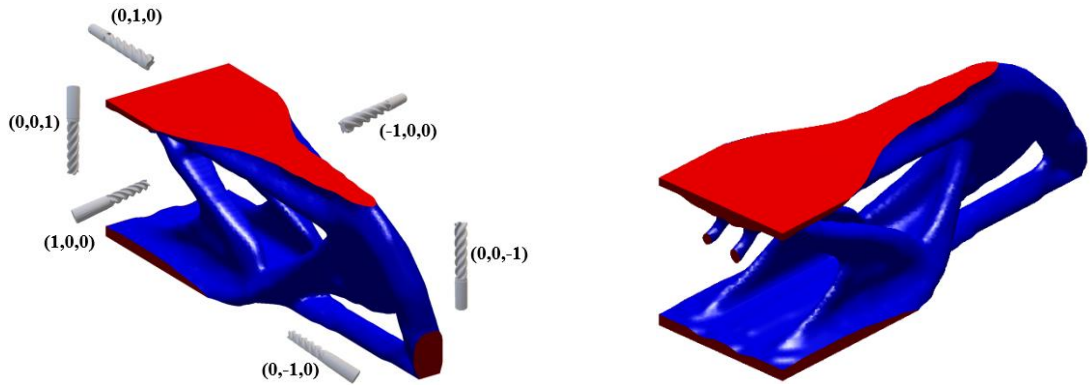

Figure 15. Multiple machining tool orientations (Compliance=11.85)

To further examine the effect of filter radius on the final optimized design, different filter sizes $r_{min}$ are selected to produce diverse designs as shown in Fig. 16. The machining tool orientations (6 directions) are shown in Fig. 15. As mentioned by Refs. [28-30], the filter size is a straightforward and effective way to control the member size of optimal design. As plotted in Fig. 16, for small radius size ($r_{min} = 1$), small and narrow hollow chambers are found, which may not be accessible by milling tool. However, the hollow chamber size increases after increasing the filter radius $r_{min}$, while the compliance values ($C$) of multiple designs are close in this case.

$(r_{min} = 1; C = 12.21)$  $(r_{min} = 2; C = 12.28)$  $(r_{min} = 3; C = 11.90)$  $(r_{min} = 5; C = 11.49)$

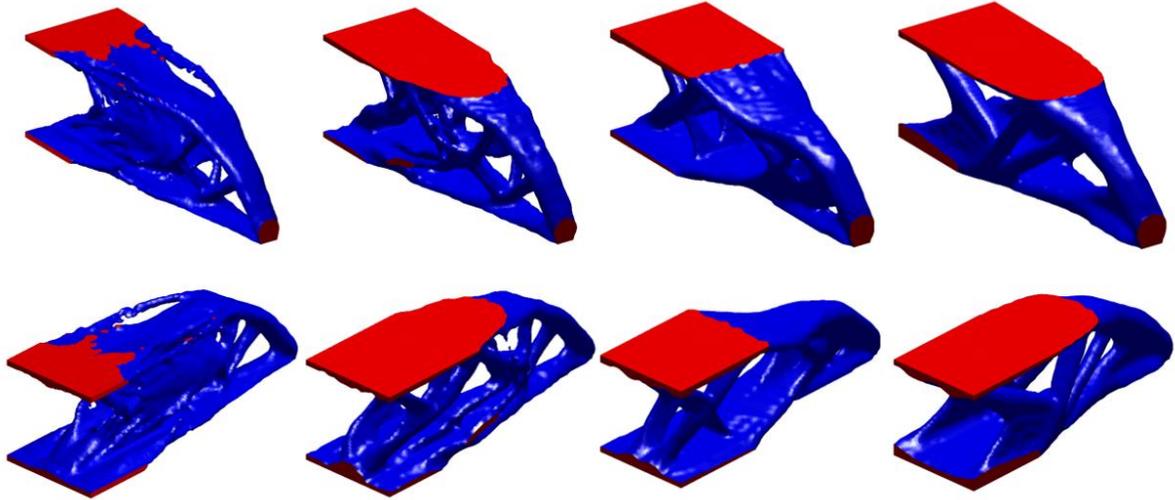

Figure 16. Machinable designs with different filter radius

### 3.2.2 Optimized Design with CES constraint

In this section, the results with CES constraint are presented. The CES angle limitation is chosen as $\theta^* = 30^o$. The filter radius is $r_{min} = 1.5$. In this example, we use the optimized design from standard density-based result to initialize the fictitous field. The initialization of fictitous domain and corresponding void field are plotted in Fig. 17. Based on numerical experiments, level set initialization method or initialization based on standard density result are better than uniform field. The problem for the uniform field initialization is that the convergence rate is very slow. Thus, we do not recommend applying uniform density distribution to initialize the design domain. Different milling directions are applied to generate multiple machinable solutions. The initial designs and optimized designs are demonstrated in Fig. 18. Obviously, compared to results in Fig. 16, no deep hole or narrow chamber is found in the final design. The optimal compliance values with CES constraint are generally higher than optimal designs without CES constraint.

(a) void field                                  (b) fictitous field

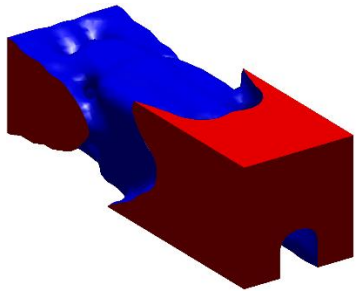 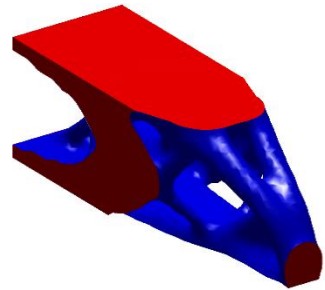

Figure 17. Design Initialization (a) void field (b) fictitous field

(a) Initial design                              (b) Optimized design (C= 23.71)

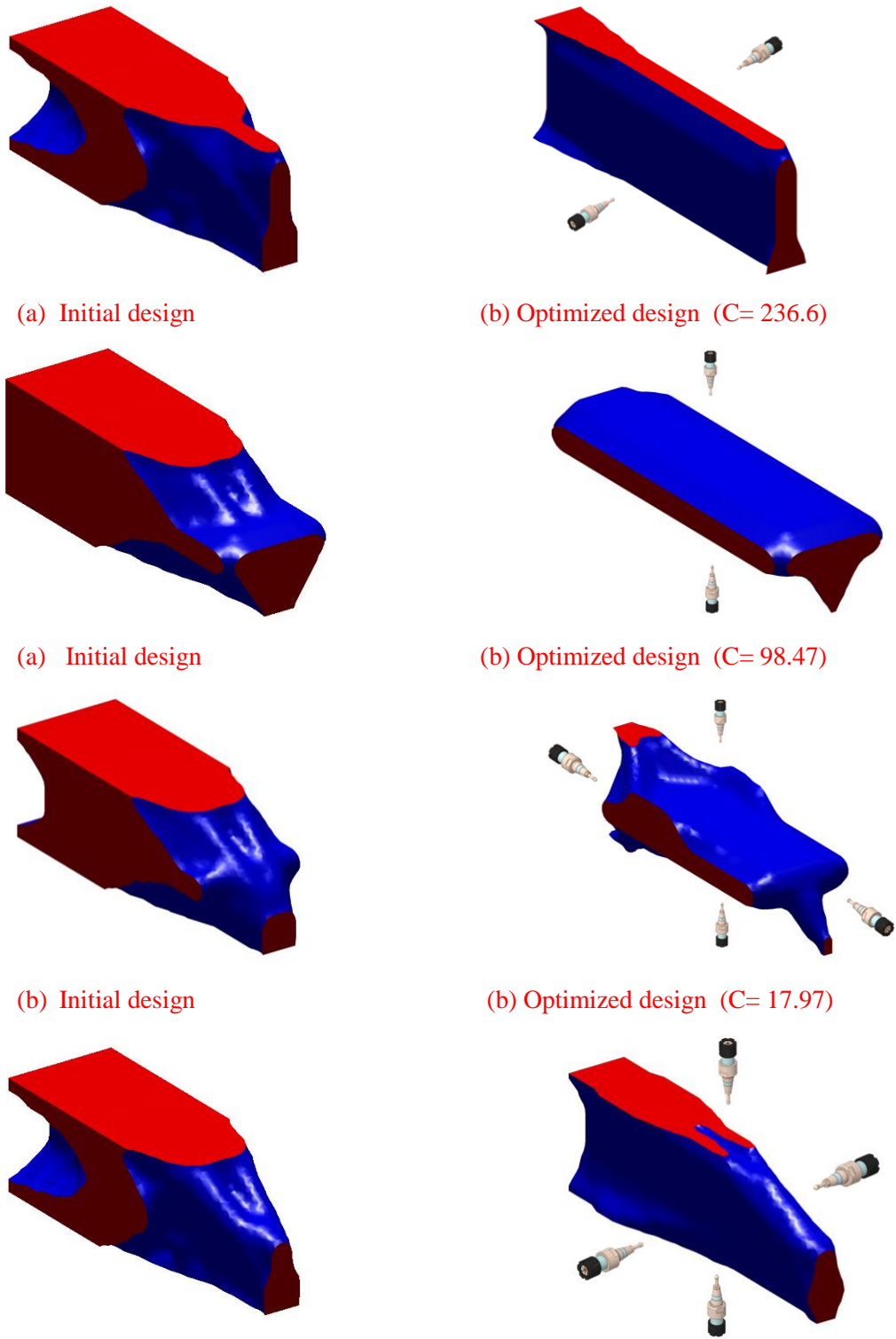

(a) Initial design     (b) Optimized design (C= 236.6)

(a) Initial design     (b) Optimized design (C= 98.47)

(b) Initial design     (b) Optimized design (C= 17.97)

Figure 18. Machinable designs with CES constraint ($\theta^* = 30^o$)

### 3.3 Machining-based Optimization for 3D MBB beam

In this section, an MBB design example is presented for compliance minimization. The unit force is applied at the center of top face. Left and right bottom edges are fully fixed as shown in Fig. 17. Due to the symmetry, only half of the MBB beam is chosen for optimization and is discretized by a $144 \times 48 \times 48$ hexahedral mesh with unit element length. The material properties are the same as the previous example. The filter size for design variables is selected as $r_{min} = 4$. The volume constraint is chosen as $\bar{V} = 0.2$. The reference result without machining constraints obtained from standard density-based method is shown in Fig. 19. This reference design is not manufacturable through machining as it contains several inaccessible internal surfaces. Similar to the previous example, we start by considering milling operation in two opposite directions. The initialization of fictitious field for optimization is demonstrated in Fig. 18. Several machinable designs are generated through the proposed method as demonstrated in Fig. 20. The machinable MBB designs have very different optimized configurations compared to the reference case, while the compliance values of machinable designs are remarkably competitive with the reference one.

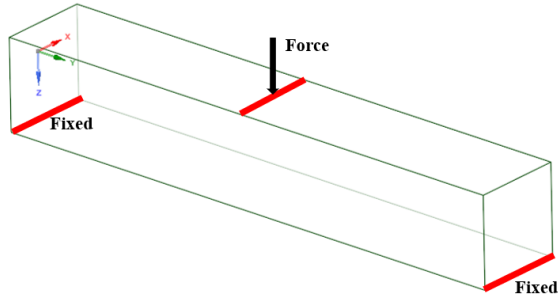

Figure 17. Machinable designs with different filter radius

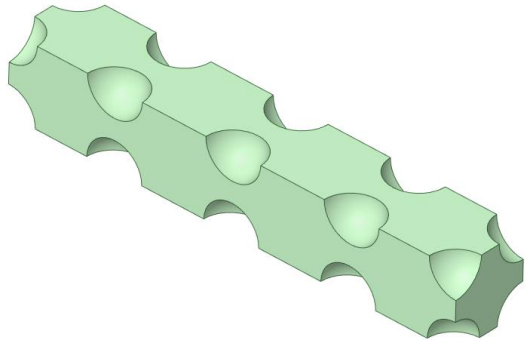

Figure 18. Initialization of fictitious field

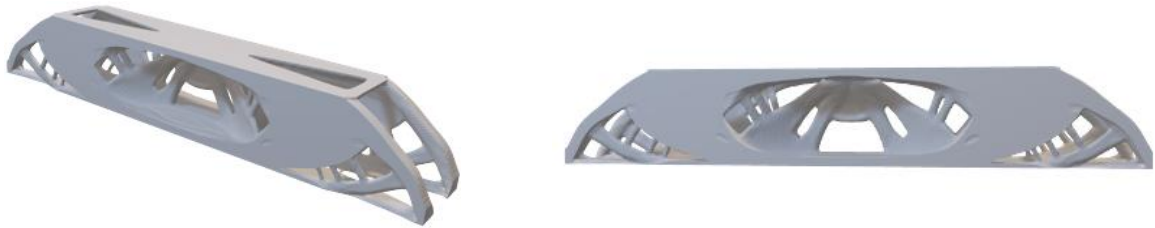

Figure 19. Reference design (Compliance=10.72)

(a) Compliance:11.38

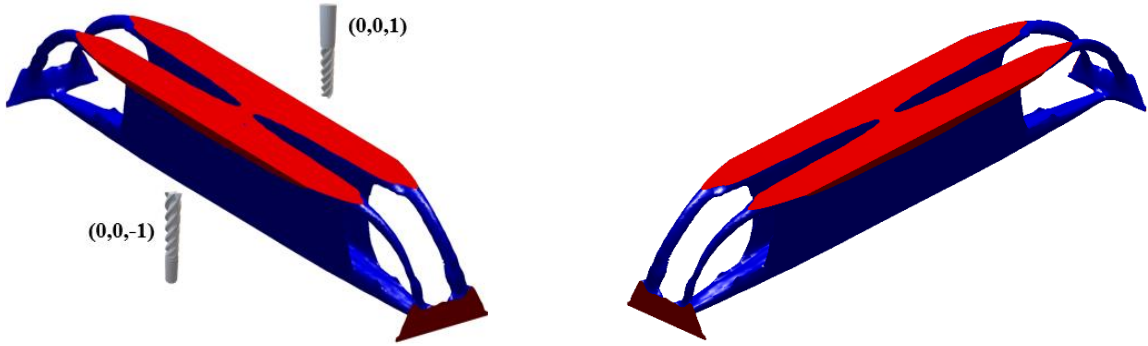

(b) Compliance: 10.21

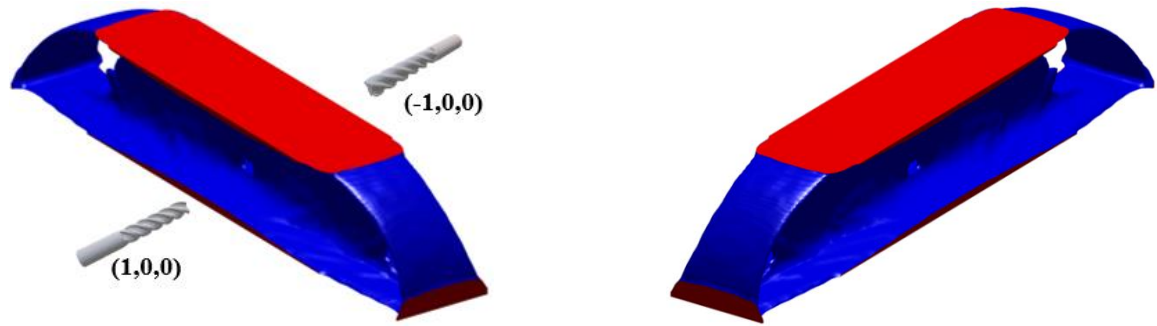

(c) Compliance: 11.09

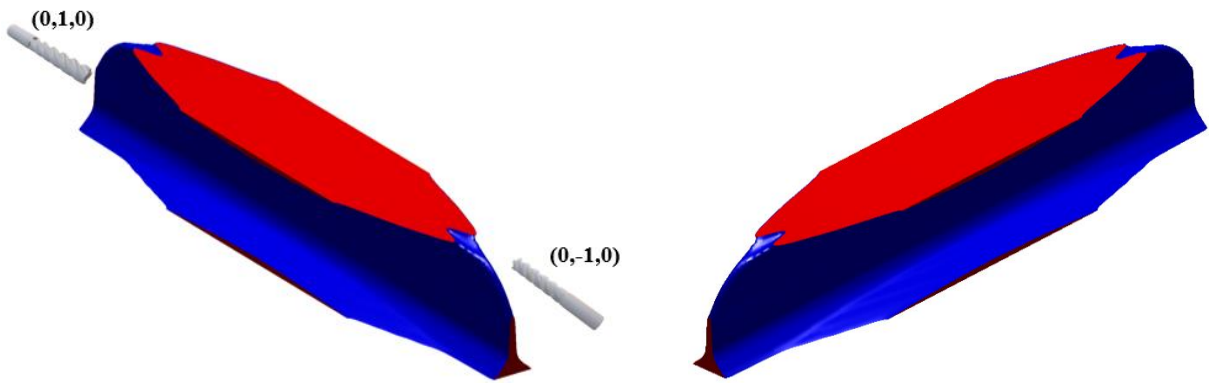

Figure 20. Two opposite milling directions

For a multi-axis milling scenario, the set of milling orientations is extended as shown in Fig. 21. As plotted in Fig. 21, the machinable design is obtained by 6 mutually orthogonal milling directions. The compliance value ($C$=10.05) of optimized design is slightly outperform than the reference design ($C$=10.72). Clearly, the reference design is only a local optimum as gradient-based method is only able to converge to a local minimum.

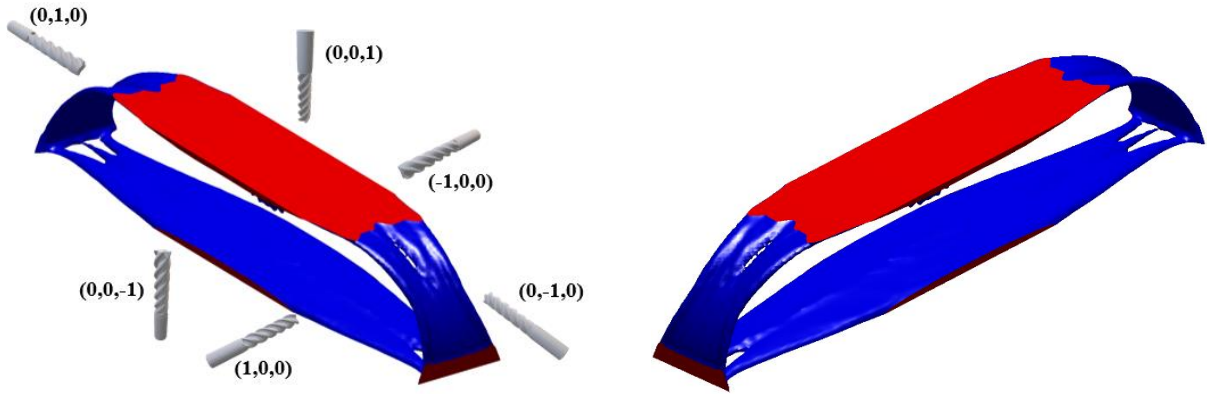

Figure 21. Multiple milling directions (Compliance: 10.05)

## 4. Conclusion

In recent years, developing advanced topology optimization methods for conventional subtractive manufacturing becomes a new trend in this field. In this paper, a novel mathematical formulation to impose multi-axis machining restrictions in the framework of density-based method is proposed. A simple density mapping method based on Heaviside function for multi-axis restrictions is demonstrated in detail, where no aggregation functions (E.g. KS-function [19]) are involved. Several 2D and 3D numerical examples are demonstrated to validate the effectiveness of proposed method. In general, some small features generated by topology optimization are difficult to reach like a thin and deep pocket. For multi-axis machining, deep and narrow areas require specialized tooling and are time-consuming to produce. These narrow areas may require additional equipment setup and increase the cost of a component. Increasing the filter radius as described in the proposed method can effectively reduce these narrow regions. Besides, the CES constraint we introduced in this paper can effectively remove these narrow or deep void features. We need to mention that the current method is assuming that the multi-axis tool's motion is only to be translational from outside the part, which is a limitation of current design method. Compared to Amir's method [16], the proposed method does not consider the fixturing devices for the machining setup, which will be further investigated and extended to realistic geometry in the future. Besides, optimizing the machining direction is another research direction which will be explored in the future.

## Acknowledgement

The financial support for this work from National Science Foundation (CMMI-1634261) is gratefully acknowledged.